\documentclass[12pt]{article}
\advance\voffset by -0.7cm
\advance\hoffset by -1.5cm
\input amssym.def
 
\def\N{{\Bbb N}}     \def\Z{{\Bbb Z}}
     \def\R{{\Bbb R}}
\def\C{{\Bbb C}}     \def\D{{\Bbb D}}
\def\H{{\Bbb H}}     

\def\semi{{\Bbb n}}
\def\SG2{{{\rm{Sp}}}(4,\Z)}
\def\SP{\mbox{{\rm{Sp}}}(2n,\Z)}
\def\SL{\mbox{SL}(2,\Z)}
\def\grad{\nabla}
\def\trace{\mbox{tr}}
\def\bin#1#2{{\left( #1\atop #2 \right)}}

\def\sq{\hbox{\rlap{$\sqcap$}$\sqcup$}}
\def\qed{\ifmmode\sq\else{\unskip\nobreak\hfil
         \penalty50\hskip1em\null\nobreak\hfil\sq
         \parfillskip=0pt\finalhyphendemerits=0\endgraf}\fi}
\def\id{{\mathchoice {\rm 1\mskip-4mu l} {\rm 1\mskip-4mu l}
        {\rm 1\mskip-4.5mu l} {\rm 1\mskip-5mu l}}}

\newtheorem{thm}{Theorem}[section]
\newtheorem{lem}[thm]{Lemma}
\newtheorem{prop}[thm]{Proposition}
\newtheorem{cor}[thm]{Corollary}
\newtheorem{df}[thm]{Definition}
\newtheorem{rem}[thm]{Remark}

\begin{document}
\date{}

\title{Rankin-Cohen Operators  
       for Jacobi and Siegel Forms}

\author{YoungJu Choie 
        \thanks{ Partially supported by KOSEF 941-01100-001-2
                 and Basic Science BSRI 96-1431.}\\
        {\small Department of Mathematics,}\\
        {\small Pohang Institute of Science \& Technology,}\\
        {\small Pohang, Korea, 790-784}\\
\and
\\
        Wolfgang Eholzer
          \thanks{Supported by the EPSRC and partially by 
                  PPARC and  EPSRC (grant GR/J73322).}
 \\
 {\small Department of Applied Mathematics and Theoretical Physics,}\\
 {\small University of Cambridge, Silver Street,}\\
 {\small Cambridge CB3 9EW, U.K.}
}

\maketitle

\abstract{
For any non-negative integer $v$ we construct explicitly 
$\lfloor \frac{v}{2}\rfloor+1$
independent covariant bilinear differential 
operators from $J_{k,m}\times J_{k',m'}$ to $J_{k+k'+v,m+m'}$.
As an application we construct a covariant bilinear
differential operator mapping $S_k^{(2)}\times S^{(2)}_{k'}$ to 
$S^{(2)}_{k+k'+v}$. Here $J_{k,m}$ denotes the space of Jacobi 
forms of weight $k$ and index $m$ and $S^{(2)}_k$ the space of 
Siegel modular forms of degree $2$ and weight $k$.
The covariant bilinear differential operators constructed are 
analogous to operators already studied in the elliptic 
case by R.\ Rankin and H.\ Cohen and we call them 
Rankin-Cohen operators.
}
\vspace{1.5cm}

\begin{center}

\end{center}

\vspace{1.0cm}
\begin{flushright}
DAMTP-96-106 \\
alg-geom/96mmnnn
\end{flushright}

\newpage

\section{Introduction and results}
\label{Introduction}

It has been known for some time how to obtain an elliptic modular form 
from the derivatives of $N$ elliptic modular forms. 
The case $N=1$ has already been studied in detail by R.\ Rankin 
in 1956 \cite{R}. For  $N=2$ H.\ Cohen has constructed 
certain covariant bilinear operators which he used to obtain 
modular forms with interesting Fourier coefficients \cite{C}.
Later, these operators were called Rankin-Cohen operators by D.\ Zagier  
who studied their algebraic relations \cite{Z}.

The main result of this paper is the explicit description of  
covariant bilinear operators for Jacobi forms and Siegel modular forms
of degree $2$.
Since they are generalisations of the Rankin-Cohen operators 
in the elliptic case we will also call them Rankin-Cohen operators.

The main theorem reads
\begin{thm}
\label{Jacobithm}
Let $f$ and $f'$ be Jacobi forms of weight and index 
$k,m$ and $k',m'$, respectively. 
For any $X\in\C$ and any non-negative integer $v$ define  
$$ [f,f']_{X,v} = \sum_{r+s+p= \lfloor v/2\rfloor \atop 
                        i+j = v-2\lfloor v/2\rfloor} 
                  C_{r,s,p}(k,k') \
                  D_{r,s,i,j}(m,m',X)\ 
                  L_{m+m'}^p\left( L_m^r(\partial_z^i f)\, 
                                   L_{m'}^s(\partial_z^j f') \right)
$$
where
\begin{eqnarray}
 D_{r,s,i,j}(m,m',X) &=& m^j (-m')^i (1+mX)^s (1-m'X)^r\, ,
\nonumber\\ 
C_{r,s,p}(k,k') &=& 
    \frac{(\alpha +r+s+p)_{s+p} }{ r! } \cdot
    \frac{(\beta  +r+s+p)_{r+p} }{ s! } \cdot
    \frac{(-(\gamma+r+s+p))_{r+s} } {p!}
\nonumber\\
&&(\alpha = k-3/2,\ \beta = k'-3/2,\ 
   \gamma = k+k'-3/2 +(v-2\lfloor v/2\rfloor) )\, ,
\nonumber
\end{eqnarray}
where $(x)_m = \prod_{0\le i \le m-1} (x-i)$, $\lfloor x\rfloor$ 
denotes the largest integer $\le x $, and, where
$L_m(f)  =  (8 \pi i m \frac{\partial}{\partial \tau}  
                    - \frac{\partial^2}{\partial z^2})f$ 
for $f$ a Jacobi form of index $m$. \\
Then $[f,f']_{X,v}$ is a Jacobi form of weight $k+k'+v$ and index
$m+m'$ and, even more, a Jacobi cusp form for $v>1$.
\end{thm}

Let us remark that some special cases of Theorem \ref{Jacobithm} 
have already been considered in the literature:
Firstly, the bilinear operator $[\cdot,\cdot]_{X,1}$,
which actually does not depend on $X$, has already been shown to map
two Jacobi forms to a Jacobi form \cite[Theorem 9.5]{EZ}.
Secondly, one bilinear operator for each even $v$ has already 
been constructed in ref.\ \cite{Ch1}: 
up to a scalar multiple this operator is equal to  
$$ \left(\frac{d}{dX}\right)^{v/2} [f,f']_{X,v}
   \qquad (v\in 2\N =\{0,2\dots \}).$$ 

For fixed $v$ and $k,m$ and $k',m'$ large enough the operators 
$[\cdot,\cdot ]_{X,v}$ ($X\in\C$) span a vector 
space of dimension $\lfloor \frac{v}{2}\rfloor+1$
This shows that the space of such Rankin-Cohen operators 
is, in general, at least $\lfloor\frac{v}{2}\rfloor+1$ dimensional.
A result of B\"ocherer \cite{Boe}, obtained by using Maa{\ss} operators, 
shows that this dimension actually equals
$\lfloor\frac{v}{2}\rfloor+1$ in general (cf.\ Theorem \ref{Boecherer}).

\medskip
One of the applications of our result is to Siegel modular forms 
of degree $2$: the bilinear operators $[\cdot,\cdot ]_{0,v}$ 
with even $v$ can be lifted to bilinear covariant differential 
operators for Siegel modular forms of degree $2$.

More precisely one has the following theorem.  
\begin{thm}
\label{Siegelthm}
Let $F$ and $F'$ be Siegel modular forms of degree $2$
and weight $k$ and $k'$, respectively.
Define, for any non-negative integer $l$,  
$$ [F,F']_{l} = \sum_{r+s+p=l} C_{r,s,p}(k,k')\
    \D^p( \D^r(F)\, \D^s(F') )
$$
with $C_{r,s,p}(k,k')$ as in Theorem \ref{Jacobithm} and
$\alpha = k-3/2$, $\beta = k'-3/2$, $\gamma = k+k'-3/2$
and where 
$\D=   \frac{4\partial^2}{\partial_{\tau_1}\partial_{\tau_2}}
      -\frac{\partial^2}{\partial_{z}^2}
$ 
with  
$Z=\left(\begin{array}{cc} \tau_1 & z \\ z &\tau_2 \end{array}\right)$ 
the variable in $\H_2$.\\
Then $[F,F']_{l}$ is a Siegel modular form of degree $2$
and weight $k+k'+2l$ and, even more, a Siegel cusp form for $l>0$.
\end{thm} 

This paper is organised as follows.
In \S\ref{Jacobiforms} we recall some results on Jacobi forms
which will be needed in the proof of Theorems \ref{Jacobithm}.
Section \ref{Siegelforms} contains the definition of Siegel modular forms 
and their relation to Jacobi forms in the degree $2$ case.
In \S\ref{Jacobiproofcom} we give a (combinatoric) proof of 
Theorem \ref{Jacobithm} and \S\ref{Jacobiproofgen}
contains a proof using generating functions.
In \S\ref{SiegelproofJacobi} we prove 
Theorem \ref{Siegelthm} using Theorem \ref{Jacobithm} and 
in \S\ref{SiegelproofTheta} we give a second, independent proof 
using theta series with spherical coefficients
and a general result of Ibukiyama (cf. Theorem \ref{Ibuthm}).
We conclude with several remarks and some open questions in 
section \S\ref{Conclusion}. 
In particular we discuss the uniqueness of the Rankin-Cohen operators 
for Siegel modular froms of degree $2$.

\section{Jacobi forms}
\label{Jacobiforms}

In this section we recall a few general results about Jacobi forms.
We first give the definition of Jacobi forms and the heat operator
(as a general reference for Jacobi forms we refer to \cite{EZ}). 

Denote by $\H$ the complex upper half plane and  define,
for holomorphic functions $f:\H \times \C \to \C$ and integers 
$k$ and $m$, 
\begin{eqnarray}
(f |_{k, m} M )(\tau,z) &=&
       (c \tau + d)^{-k} e^{2 \pi i m (\frac{-c z^2}{c \tau + d})}
       f(\frac{a \tau + b}{c \tau +d}, \frac{z}{c \tau +d}),
\nonumber\\
(f|_{m} Y) (\tau,z) &=& 
      e^{2 \pi i m (\lambda^{2} \tau + 2\lambda z)}
      f(\tau, z+\lambda \tau + \nu)
\nonumber
\end{eqnarray}
where $\tau\in\H$, $z\in\C$,  
$M = \left(\begin{array}{cc} a & b\\ c & d \end{array} \right) \in \SL$
and $Y=(\lambda,\nu)\in\Z^2$.

Using these slash actions the definition of Jacobi forms is as follows.
\begin{df}
A Jacobi form of weight $k$ and index $m$ ($k, m \in \N$)
is a holomorphic function  $f:\H\times\C\to\C$ satisfying 
$$ (f |_{k, m} M)(\tau, z) = f(\tau,z), \qquad
   (f|_m  Y) (\tau,z)= f(\tau,z)
$$ 
for all $M\in\SL$ and $Y\in\Z^2$ and such that it has a Fourier 
expansion of the form
$$ f(\tau, z)  =    
                 \sum_{ \begin{array}{cc}
                        n = 0 \\
                        r \in \Z, r^{2} \leq 4nm
                        \end{array}
                      }^{\infty}
                   c(n, r)  q^n\zeta^r,
$$
where $q = e^{2 \pi i \tau}$ and $\zeta = e^{2 \pi i z}$.
If $f$ has a Fourier expansion  of the same form but with 
$r^2 < 4nm$ then $f$  is called  a Jacobi cusp form 
of weight $k$ and index $m$.  
\end{df}
We denote by $J_{k,m}$ the (finite dimensional) vector space of 
all Jacobi forms of weight $k$ and index $m$ and by $J_{k,m}^{cusp}$ 
the vector space of all Jacobi cusp forms of weight $k$ 
and index $m$.

Our main result (Theorem \ref{Jacobithm}) involves the heat operator 
which has already been studied in \cite{EZ} to connect Jacobi forms 
and elliptic modular forms and in ref.\ \cite{Ch1,Ch2} in the context
of bilinear differential operators ({\it cf.} the remark after Theorem
\ref{Jacobithm}).  
\begin{df}
For any non-negative integer $m$ the heat operator
$L_m$ is defined by 
$$ L_m(f)  =  \left(8 \pi i m \frac{\partial}{\partial \tau}  
                    - \frac{\partial^2}{\partial z^2}
              \right)(f) 
   \qquad(f\in J_{k,m}).
$$
\end{df}

Finally, let us mention a result of B\"ocherer \cite{Boe}.
\begin{thm}
\label{Boecherer}
For fixed $v$ and $k,m,k',m'$ large enough the vector space of all 
covariant bilinear differential operators mapping 
$J_{k,m}\times J_{k',m'}$ to $J_{k+k'+v,m+m'}$ has  dimension 
$\lfloor v/2\rfloor+1$.
\end{thm}
Note that Theorem \ref{Jacobithm} describes a basis of this space 
explicitly.

\section{Siegel modular forms} 
\label{Siegelforms}

In this section we recall a few basic facts about Siegel modular forms
and, in particular, the construction of Siegel modular forms using 
theta series with spherical coefficients (Theorem \ref{thetathm}).
Furthermore, we describe the connection between Siegel modular forms of
degree $2$ and Jacobi forms (Theorem \ref{Siegelexp}).
Finally, we mention (a special case of) a result of Ibukiyama  
(Theorem \ref{Ibuthm}) which will be needed in the proof of Theorem
\ref{Siegelthm} in section \S\ref{SiegelproofTheta}.
The reader may take ref.\ \cite{Frei} as a general reference for 
Siegel modular forms.

For any holomorphic function $f$ on the Siegel upper half plane
$\H_n$, {\it i.e.}\ the space of complex
symmetric $n\times n$ matrices with positive definite imaginary part,
and $M\in \SP$ define
$$ (f|_M^k)(Z) = f(MZ) \det(CZ+D)^{-k} $$
where $Z\in\H_n$, 
$M=\left( \begin{array}{cc} A & B \\ C & D \end{array}\right)$
with $n\times n$ matrices $A,B,C,D$, and, where
$MZ = (AZ+B)(CZ+D)^{-1}$.

Then the definition of Siegel modular forms is given as follows.
\begin{df}
Let $n$ be a positive integer greater than $1$.
Then a holomorphic function $f$ on $\H_n$ is called a Siegel modular 
form of degree $n$ and weight $k$  if 
$$ (f|_M^k)(Z) = f(Z)\qquad {\rm{for}}\ {\rm{all}}\  M\in\SP. $$
\end{df}
We denote the space of all Siegel modular forms of degree $n$ and
weight $k$ by $S^{(n)}_k$. 

Note that for $n=1$ one has to add a further condition on $f$ 
in order to obtain the usual definition of modular forms.

The connection between Jacobi forms and Siegel modular forms of 
degree two becomes clear by the following theorem (see {\it e.g.} 
\cite[Theorem 6.1]{EZ}).
\begin{thm}
\label{Siegelexp}
Let $F$ be a Siegel modular for degree $2$ and weight $k$ and
write the Fourier development of $F$ with respect to $\tau_2$
in the form
$$ F(Z) = F(\tau_1,z,\tau_2) = \sum_{m=0}^\infty 
                         f_m(\tau_1,z) e^{2\pi i m\tau_2}
$$
with $Z =\left( \begin{array}{cc} \tau_1& z \\ z & \tau_2 \end{array}\right)$
the variable in $\H_2$.\\ 
Then, for each non-negative integer $m$, the function $f_m$ is a 
Jacobi form of weight $k$ and index $m$.
\end{thm}

\begin{rem}
\label{Siegelcusp}
Note that for general degree $n$ a Siegel cusp form is 
a Siegel modular form which is contained in the kernel of 
the Siegel operator $\Phi$ (for more details see {\it e.g.}
\cite{Frei}). In the case of degree $2$ a Siegel cusp form $F$ 
is a Siegel modular form whose Jacobi-Fourier expansion is of the form 
$F= \sum_{m>0} f_m(\tau_1,z) e^{2\pi i m\tau_2}$, {\it i.e.}\ 
the first coefficient in the Jacobi-Fourier expansion is identically 
zero.
\end{rem}

To recall some facts about theta series with spherical coefficients 
we introduce the notion of spherical polynomials first.
\begin{df}
A spherical polynomial $P$ of weight $w$ in a matrix variable 
$X\in M_{m,n}$ is a polynomial (in the matrix elements of 
$X$, {\it i.e.}\ a polynomial in  
$\C[x_{ij}]_{1\le i \le m \atop 1\le j \le n}$)  satisfying \\
1) $P(X A) = \det(A)^w\, P(X)$ for all $A\in M_{n,n}$, \\
2) $\Delta P = \sum_{i,j} \frac{\partial^2}{(\partial x_{ij})^2} P =0$.
\end{df}

Then one has (see {\it e.g.} \cite[p.\ 161]{Frei})
\begin{thm}
\label{thetathm}
Let $P$ be a spherical polynomial of weight $k$ in a matrix variable 
$X\in M_{m,n}$ and let $S\in M_{m,m}(\Z)$ be a symmetric, positive, even 
and unimodular matrix.
Then the function 
$$ \theta_{S,P}(Z) = \sum_{G\in M_{m,n}(\Z)} 
                       P(S^{1/2}G)\, e^{\pi i \, \trace(G^t S G Z)}
$$
is a Siegel modular form of degree $n$ and weight $m/2+k$.
\end{thm}

Finally, let us mention (a special case of) a result of Ibukiyama 
(Corollary 2 (2) of ref.\ \cite{Ibu} with $r=2$) which we use  
to prove Theorem \ref{Siegelthm} in \S\ref{SiegelproofTheta}.
\begin{thm}
\label{Ibuthm}
Let $P$ be a spherical polynomial of even weight $d$ in the matrix 
variable $(X,X')^t \in M_{m+m',n}$ which can be written as 
$P(X,X') = \tilde Q(X^t X, {X'}^t X')$ for some polynomial $\tilde Q$.
Set ${\cal D} = \tilde Q(\partial_{\nu,\mu},\partial_{\nu',\mu'}')$
where 
$\partial_{\nu,\mu} = (1+\delta_{\nu,\mu})
 \frac{\partial}{\partial z_{\nu,\mu}}$
and 
$\partial_{\nu',\mu'}' = (1+\delta_{\nu',\mu'})
 \frac{\partial}{\partial z'_{\nu',\mu'}}$
with $Z = (z_{\nu,\mu})$
and  $Z' = (z'_{\nu',\mu'})$ for 
$1\le \nu,\mu,\nu',\mu' \le n$.\\
Then, for any two Siegel modular forms $F(Z)$ and $F'(Z')$ 
of degree $n$ and weight $k$ and $k'$, respectively, the function
${\cal D}(F(Z) F'(Z'))|_{Z=Z'}$ is a Siegel modular form  of 
degree $n$ and weight $k+k'+d$.
\end{thm}
Note that Theorem \ref{Ibuthm} essentially says that 
if a bilinear differential operator maps all pairs of theta 
series to theta series with spherical coefficients then 
it even maps all pairs of Siegel modular forms to Siegel 
modular forms ({\it cf.} the discussion at the end of \S6).

\section{A combinatorial proof of Theorem \ref{Jacobithm}}
\label{Jacobiproofcom}

Let us now give the proof of our main theorem which will, in particular,
imply Theorem \ref{Siegelthm} (see \S\ref{SiegelproofJacobi} for details).

To prepare the proof of Theorem \ref{Jacobithm} we need the 
following three lemmata. 
\begin{lem}
\label{Lmdis}
Let $f$ be a holomorphic function on $\H$ and $g$ a holomorphic
function on $\H\times\C$.
Then, for each non-negative integer $r$,  one has 
$$ L_m^r(fg) = \sum_{j=0}^r(8\pi i m)^{r-j} \bin{r}{j}\,
                           (\partial_{\tau}^{r-j} f)\,
                           (L_m^j g)
$$
where $\tau$ is the variable in $\H$.
\end{lem}
{\it Proof.}
We prove the formula by induction. 
Firstly,  note that for $r=1$ one has
$$ L_m(fg) = (8\pi i m \partial_{\tau}-\partial_{z}^2)(fg)
          =  8\pi i m (\partial_{\tau}f)g +  
             8\pi i m f (\partial_{\tau}g) - 
             f\partial_{z}^2 g
          = 8\pi i m (\partial_{\tau}f) g +  f L_m g.
$$
Secondly, assume that the formula is valid for some $r$. 
Then we find
\begin{eqnarray}
 L_m^{r+1}(fg) &=& L_m\left( \sum_{j=0}^r (8\pi i m)^{r-j} \bin{r}{j}\, 
                           (\partial_{\tau}^{r-j} f)\,
                           (L_m^j g)\right)
\nonumber\\
              &=& \sum_{j=0}^r \Big( 
                           (8\pi i m)^{r-j+1}\bin{r}{j}\, 
                           (\partial_{\tau}^{r-j+1} f)\,
                           (L_m^j g)
\nonumber\\
&&                           \qquad+
                           (8\pi i m)^{r-j}\bin{r}{j}\, 
                           (\partial_{\tau}^{r-j} f)\,
                           (L_m^{j+1} g)
                           \Big)
\nonumber\\
              &=& \sum_{j=0}^{r+1}
                            (8\pi i m)^{r+1-j}\bin{r+1}{j}\, 
                           (\partial_{\tau}^{r+1-j} f)\,
                           (L_m^j g)
\nonumber
\end{eqnarray}
and the lemma becomes obvious.
\qed

The second lemma we will need is Lemma 3.1 of ref.\ \cite{Ch1} (note 
that the normalisation of $L_m$ in {\it loc. cit.} differs from ours 
by a factor of $(2\pi i )^2$).
\begin{lem}
\label{Lmcom}
Let $f$ be a holomorphic function on $\H\times\C$.
Then, for any non-negative integer $r$, one has
$$
L_m^r(f)|_{k+2r,m}M = \sum_{j=0}^r (8\pi i m)^{r-j}
                                   \frac{r!}{j!} 
                                   \bin{k-3/2+r}{r-j}\,
                                \left(\frac{c}{c\tau+d}\right)^{r-j}\,
                                L_m^j(f|_{k,m} M)
$$
for all 
$M=\left(\begin{array}{cc} a&b\\ c&d \end{array}\right)\in{\mbox{SL}}(2,\Z)$.\\
Furthermore, one has 
$$ L_m(f)|_m Y = L_m(f|_m Y) $$
for all $Y\in \Z^2$.
\end{lem}
{\it Proof.} 
The proof of the first formula is a simple exercise and can be found 
in \cite{Ch1}. The second formula becomes obvious after a short 
calculation using only the definition of the $|_m Y$ action
\qed

The third lemma is essentially equivalent to Theorem 9.5 of ref.\ \cite{EZ}.
\begin{lem}
\label{oddlem}
Let $f$ and $f'$ be holomorphic functions on $\H\times\C$.
Then, for $z\in \C$,  one has
\begin{eqnarray}
 &&\left( m' (\partial_z f)\, f' -
        m  f\, (\partial_z f')
 \right)|_{k+k'+1,m+m'} M = 
\nonumber\\
&&\qquad\qquad\qquad \qquad
         m' (\partial_z (f|_{k,m} M))\, (f'|_{k',m'} M) -
         m  (f|_{k,m} M)\, (\partial_z (f'|_{k',m'} M))\, ,
\nonumber \\ 
&&\left( m' (\partial_z f)\, f' -
        m  f\, (\partial_z f')
 \right)|_{m+m'} Y =
\nonumber \\
&&\qquad\qquad\qquad\qquad
        m' (\partial_z (f|_{m} Y))\, (f'|_{m'} Y) -
        m  (f|_m Y)\, (\partial_z (f'|_{m'} Y))
\nonumber 
\end{eqnarray}
for all $M\in\SL$ and all $Y\in\Z^2$.
\end{lem}
{\it Proof.} The two formulas can be obtained by a straightforward 
calculation. \qed
\bigskip

We are now ready to prove our main theorem.\\
{\it Proof of Theorem \ref{Jacobithm}.}
The definition of Jacobi forms contains essentially two 
parts: the invariance under the slash actions and 
an expansion condition. 
If $f$ and $f'$ are Jacobi forms then a simple computation shows that 
the latter condition is always satisfied for $[f,f']_{X,v}$ and,
even more, that $[f,f']_{X,v}$ satisfies the expansion condition of a 
Jacobi cusp form for $v>1$.
Therefore, it only remains to check the invariance under 
the slash actions.  

We consider first the case of even $v=2v'$.
Since, by Lemma  \ref{Lmcom},  the $|_m Y$ action commutes with $L_m$ 
we only have to show that 
$$ ([f,g]_{X,v})|_{k+k'+v,m+m'} M = [f|_{k,m} M,g|_{k',m'} M]_{X,v}$$
for all $M\in \SL$.
Using Lemma \ref{Lmdis} and \ref{Lmcom} to calculate the left hand side 
we obtain that this equation is equivalent to 
\begin{eqnarray}
&&\sum_{r+s+p=v'} C_{r,s,p}  \frac{r!s!p!}{j!j'!j''!} 
                            \bin{\alpha+r}{r-j'}
                            \bin{\beta +s}{s-j''}
                            \bin{\gamma+v'+j'+j''}{p-j}
\nonumber\\
&& \qquad \times \quad
                            (1+mX)^s (1-m'X)^r 
                            m^{r-j'}m'^{s+p-j-j''}
                            (1+m/m')^{p-j} 
\nonumber\\
&& = \delta_{v',j+j'+j''}\, C_{j',j'',j} (1+mX)^{j''} (1-m'X)^{j'}.
\nonumber
\end{eqnarray}
Some simple manipulations show that this equation is equivalent to
\begin{eqnarray}
&& \sum_{r+s+p=v'} \frac{(v'-j')!(v'-j'')!(v'-p)!}{(r-j')!(s-j'')!(v'-j)!}
                  \bin{-(\gamma+v')}{v'-p}
                  \bin{\gamma+v'+j'+j''}{p-j}
\nonumber\\
&&\qquad \times \quad 
                  m^{r-j'} m'^{s+p-j-j''}(1+m/m')^{p-j}
                 (1+mX)^{s-j'} (1-m'X)^{r-j''}
\nonumber\\
&&  = \delta_{v',j+j'+j''}\, \bin{-(\gamma+v')}{v'-j}.
\nonumber
\end{eqnarray}
To show this equality we view it as an equation 
between polynomials in $\gamma$ (of degree $v'-j$).
It is easy to see that both sides coincide for 
$\gamma+v'+j'+j''=0$. Therefore, it is enough to show that 
both sides agree for all $\gamma =  -v'-x$ $(0\le x <v'-j$).
Note that for these values of $\gamma$ the right hand side 
obviously vanishes.
With $A = \frac{m(1-m'X)}{m'(1+mX)}$ the left hand side 
becomes 
\begin{eqnarray}
 && c\  \sum_{r,p} \bin{v'-p-j'-j''}{r-j'}
             \bin{x-j'-j''}{v'-p-j'-j''}
             \bin{j'+j''-x}{p-j}
                  A^r (A+1)^{p-j}
\nonumber\\
&& \qquad     = c'\ \sum_{i'}\sum_{r,p} 
                     \bin{j+j'+j''-v'}{x+p-v'}
                     \bin{v'-j-j''-i'}{v'-r-p-j''}
                     \bin{i'-j'}{r-j'}
                     A^{i'}
\nonumber
\end{eqnarray}
for some (non-zero) factors $c$ and $c'$.
The product of the three binomial coefficients is the factor
in front of 
$$ Z^{x+p-v'}\cdot Z^{v'-r-p-j''} \cdot Z^{r-j'} = Z^{x-j'-j''}$$
in 
$$(1+Z)^{j+j'+j''-v'} \cdot (1+Z)^{v'-j-j''-i'} \cdot 
  (1+Z)^{i'-j'} = 1.$$
Hence we find that our expression is zero unless $x=j'+j''$. 
However, in this case our expression can only be non-zero if
$v'=j+j'+j''$ which is not allowed since $x=j'+j''$ has to be 
strictly less than $v'-j=j'+j'' = x$.
Hence also the left hand side is equal to zero for 
$\gamma =  -v'-x$ $(0\le x <v'-j$) so that we have proved the 
desired equality. 
This proves the theorem for even $v=2v'$.

For odd $v=2v'+1$ note that the summand for fixed $r,s,p$ is equal to 
$$ C_{r,s,p} (1+mX)^s (1-m'X)^r L_{m+m'}^p
   \left( m' (\partial_z L_m^r(f))\,  L_{m'}^s(f') -
          m L_m^r(f)\, (\partial_z L_{m'}^s(f')) \right).
$$
The expression inside the $L_{m+m'}^p$ is exactly of the form 
considered in Lemma \ref{oddlem} so that we obtain 
$$ ([f,f']_{X,v})|_{m+m'} Y = [f|_{m} Y,f'|_{m'} Y]_{X,v}.$$
Therefore, it only remains to show the invariance with 
respect to the other slash action as in the even case.
Note, however, that using Lemma \ref{oddlem} we have to check 
exactly the same combinatorial identity as in the case of 
even $v$  but with $\gamma = k+k'-1/2$ instead of 
$\gamma = k+k'-3/2$. This completes the proof of Theorem \ref{Jacobithm}.
\qed

\section{A proof of Theorem \ref{Jacobithm} using generating functions}
\label{Jacobiproofgen}

In this section we give a second, independent proof of Theorem 
\ref{Jacobithm}. Instead of proving the theorem directly we use 
some results of ref.\ \cite{Ch2} (Theorem 3.1 and Corollary 3.1) 
on generating functions. Let us first recall these results.
\begin{thm}
\label{Choiethm}
Let $\tilde f(\tau,z;W)$ be a formal power series in $W$, {\it i.e.} 
$\tilde f$ can be written as
$\tilde f(\tau,z;W) = \sum_{\nu=0}^\infty \chi_\nu(\tau,z)\, W^\nu$,
satisfying the functional equation 
$$ \tilde f(\frac{a\tau+b}{c\tau+d},\frac{z}{c\tau+d};\frac{W}{(c\tau+d)^2}) =
   (c\tau+d)^K e^{2\pi i M \frac{cz^2}{c\tau+d}}
               e^{8\pi i M\frac{cW}{c\tau+d}}\ \tilde f(\tau,z;W)
$$
for some integers $K$ and $M$ and 
all $\left(\begin{array}{cc} a &b\\ c& c\end{array}\right)\in
\SL$.
Furthermore, assume that the coefficients $\chi_\nu$ are holomorphic 
functions  on $\H\times\C$ with a Fourier expansion of the form
$$ \chi_\nu(\tau,z) = \sum_{r,n\in\Z\atop r^2\le 4mn}
                      c(n,r) q^n \xi^r
   \qquad( q = e^{2\pi i \tau}, \xi = e^{2\pi i z})
$$
satisfying $ \chi_\nu|_m Y = \chi_\nu$ for all $Y\in \Z^2$.\\
Then, for each non-negative integer $\nu$,  the function $\zeta_\nu$ defined by
$$ \zeta_\nu(\tau,z) = \sum_{j=0}^\nu  
                      \frac{(-(K-3/2+\nu))_{\nu-j}}{j!}\ 
                      L^{j}_M(\chi_{\nu-j})
$$
is a Jacobi form of weight $K+2\nu$ and 
index $M$.
\end{thm}
An immediate consequence of this theorem is the following corollary 
\cite[Corollary 3.1]{Ch2}.
\begin{cor}
\label{Choiecor}
Let $f(\tau,z)$ be a Jacobi form of weight $k$ and index $m$.
Then 
$$ \tilde f(\tau,z;W) = \sum_{\nu=0}^\infty
      \frac{1}{\nu!(K-3/2+\nu)!}\  L^\nu_M(f)\ W^\nu
$$ 
satisfies the functional equation
$$ \tilde f(\frac{a\tau+b}{c\tau+d},\frac{z}{c\tau+d};\frac{W}{(c\tau+d)^2}) =
   (c\tau+d)^k e^{2\pi i m \frac{cz^2}{c\tau+d}}
               e^{8\pi i m\frac{cW}{c\tau+d}}\ \tilde f(\tau,z;W)
$$
for all $\left(\begin{array}{cc} a &b\\ c& c\end{array}\right)\in\SL$.
\end{cor}

Using these two results we are now able to prove Theorem \ref{Jacobithm}.\\
{\it Proof of Theorem \ref{Jacobithm}}.
Let $f$ and $f'$ be Jacobi forms of weight and index $k,m$ and 
$k',m'$, respectively.
Denote by $\tilde f(\tau,z;W)$ and $\tilde f'(\tau,z;W)$ the formal 
power series associated to $f$ and $f'$ as in Corollary
\ref{Choiecor}, respectively.
Then, for any fixed complex number $X$, the formal power series 
$\tilde F_X(\tau,z;W)$ defined by 
$$\tilde F_X(\tau,z;W) =  \tilde f (\tau,z;(1+m'X) W)\
                          \tilde f'(\tau,z;(1-m X) W) $$ 
satisfies, by Corollary \ref{Choiecor}, the functional equation 
stated in Theorem \ref{Choiethm} with $K=k+k'$ and $M=m+m'$. 
Furthermore, it is simple to check that its coefficients satisfy 
the expansion condition assumed in Theorem \ref{Choiethm} and, 
by Lemma \ref{Lmcom}, are invariant under the $|_m Y$ action.
Hence the function $\tilde F_X(\tau,z;W)$ satisfies all assumptions of 
Theorem \ref{Choiethm} so that the corresponding functions 
$\zeta_\nu(\tau,z)$ are Jacobi forms of weight $k+k'+2\nu$ and index 
$m+m'$. It is now a simple  exercise to see that these functions are just 
constant multiples of the Rankin-Cohen operators $[f,f']_{X,2\nu}$. 
This proves Theorem \ref{Jacobithm} for even $v=2\nu$.

For the case of odd $v=2\nu+1$ consider
the function $\tilde G_X(\tau,z;W)$ defined by
\begin{eqnarray}
\tilde G_X(\tau,z;W)
   &=& \  m'(\partial_z \tilde f (\tau,z;(1+m'X) W))\
        (\tilde f'(\tau,z;(1-m X) W))
\nonumber\\
   &&  -m (\tilde f (\tau,z;(1+m'X) W))\
        (\partial_z \tilde f'(\tau,z;(1-m X) W))
\nonumber
\end{eqnarray}
where $X$ is a fixed complex number.
Using again Lemma \ref{Lmcom} and Corollary \ref{Choiecor} 
we find that the function  $\tilde G_X(\tau,z;W)$ satisfies
the functional equation of Theorem \ref{Choiethm}
with $K=k+k'+1$ and $M=m+m'$. 
By the same calculation as in the case of even $v$ the coefficients 
of $\tilde G_X(\tau,z;W)$ satisfy the expansion condition and are, by
Lemma \ref{Lmcom}, invariant under the $|_{m+m'} Y$ action.
Therefore, the corresponding functions $\zeta_\nu(\tau,z)$ are Jacobi 
forms of weight $k+k'+2\nu+1$ and index $m+m'$. These functions 
are just constant multiples of $[f,f']_{X,2\nu+1}$.
This completes the proof of Theorem \ref{Jacobithm}.
\qed

\section{A proof of Theorem \ref{Siegelthm} using Theorem \ref{Jacobithm}}
\label{SiegelproofJacobi}

In this section we give a proof of Theorem \ref{Siegelthm} 
using Theorem \ref{Jacobithm}.

\smallskip\noindent
{\it Proof of Theorem \ref{Siegelthm}.}
First we recall some well known facts about Siegel modular forms
of degree $2$ which we need in the proof. Using  
$$Z = \left(\begin{array}{cc} \tau_1 & z\\ z & \tau_2 \end{array}\right)$$
for a variable in $\H_2$ 
the `$|^k_M$' action for the whole group $\SG2$ 
(introduced in \S3) is generated by the following three 
slash actions (see {\it e.g.} page 73 of \cite{EZ})
\begin{eqnarray}
f|^k_A(\tau_1,z,\tau_2) &=& (c\tau_1+d)^{-k} 
     f(\frac{a\tau_1+b}{c\tau_1+d},\frac{z}{c\tau_1+d}, 
       \tau_2-\frac{cz^2}{c\tau_1+d})
\nonumber \\
f|_{(\lambda,\nu)}(\tau_1,z,\tau_2) &=& 
     f(\tau_1,z+\lambda \tau_1+\nu,\tau_2+2\lambda z +\lambda^2\tau_1)
\nonumber \\
f|_{t}(\tau_1,z,\tau_2) &=& f(\tau_2,z,\tau_1)
\nonumber
\end{eqnarray}
where $A = \left(\begin{array}{cc} a&b\\ c&d \end{array}\right)$ is in
${\mbox{SL}}(2,\Z)$ and $\lambda,\nu\in\Z$.  
Here we have identified, for any function $f$ on $\H_2$, 
$f(\tau_1,z,\tau_2)$ with $f(Z)$. 

For the proof of the theorem note that a simple calculation, using only
the definition of the slash actions, shows that
the operator $\D$ commutes with the $|_{(\lambda,\nu)}$ and $|_{t}$ action,
respectively. 
Secondly, it is essential to realise that 
for a Jacobi form $f$ of index $m$
$$ \D( f \tilde q^m) = L_m(f) \tilde q^m $$
where $\tilde q =e^{2\pi i \tau_2}$
and  $f$ is considered to be a function of $\tau_1$ and $z$.

Consider now two Siegel modular forms $F$ and $F'$ 
of weight $k$ and $k'$, respectively.
By Theorem \ref{Siegelexp} we can write $F$ and $F'$ as 
$$ F  = \sum_{m \ge 0} f_m     \tilde q^m,\qquad
   F' = \sum_{m'\ge 0} f'_{m'} \tilde q^{m'}
$$
where $f_m$ and $f'_{m'}$ are Jacobi forms of weight $k$ and $k'$ 
and index $m$ and $m'$, respectively.
Then
$$[F,F']_l = \sum_{m,m'\ge 0} [f_m,f_{m'}]_{X=0,2l}\  \tilde q^{m+m'}.$$

By Theorem \ref{Jacobithm} the functions $[f_m,f_{m'}]_{X=0,v}$
are Jacobi forms of weight $k+k'+2l$ and index $m+m'$ so that 
the right hand side looks like the expansion of a Siegel form
of degree $2$ and of weight $k+k'+2l$. In order to show that 
this is indeed the case the only additional property to be checked 
is that $[F,F']_l$ is symmetric in $\tau_1$ and $\tau_2$.
This, however, is obvious from the definition of  $[F,F']_l$ and the
fact that $F$ and $F'$ are themselves Siegel modular forms.
Hence, we have proved that $[\cdot,\cdot]_l$ maps 
$S^{(2)}_k\times S^{(2)}_{k'}$
to $S^{(2)}_{k+k'+2l}$.
Finally, note that the first coefficient in the
Jacobi-Fourier expansion of $[F,F']_l$ is $[f_0,f_{0'}]_{X=0,2l}$
which is identically zero for $l>0$.
This implies that the image of $[\cdot,\cdot]_l$ is contained in the 
space of Siegel cusp forms for $l>0$ ({\it cf.} Remark \ref{Siegelcusp}).
 \qed

\section{A proof of Theorem \ref{Siegelthm} using theta series}
\label{SiegelproofTheta}

In this section we give a proof of Theorem \ref{Siegelthm} 
using theta series with spherical coefficients and Theorem
\ref{Ibuthm}.
This proof also shows that the Rankin-Cohen operator $[\cdot,\cdot]_l$ 
is (up to multiplication by a constant) the only covariant bilinear 
differential operator for Siegel modular forms of degree $2$ which can
be written in terms of the differential operator $\D$.
The proof does not use Theorem \ref{Jacobithm} 
but instead Theorem \ref{Ibuthm} which essentially allows to assume 
that $F$ and $F'$ both can be written as theta series with 
harmonic coefficients.
The proof is of independent interest since part of the calculations 
are valid for arbitrary degree $n$ of the Siegel modular forms 
involved.

Firstly, we define the differential operator $\D$ for 
general degree $n$  
$$ \D = \det(\partial_{\nu,\mu}) $$
where 
$\partial_{\nu,\mu} = (1+\delta_{\nu,\mu})
 \frac{\partial}{\partial z_{\nu,\mu}}
$
with $Z = (z_{\nu,\mu})$ the usual variable in the Siegel 
half plane $\H_n$.
This operator has already been considered in 
the context of Siegel modular forms (see {\it e.g.} \cite{Frei}).

Secondly, we study the action of $\D$ on theta series.
\begin{lem}
\label{trilem}
Let $A$ be a symmetric matrix in $M_{n,n}$ and $Z$ a variable in $\H_n$.
Then one has
$$ \D\left( e^{\pi i\,  \trace(A Z)} \right) =
   (2\pi i)^n \, \det(A) \ e^{\pi i\,  \trace(A Z)}
$$  
\end{lem}
{\it Proof.} The equality follows directly from 
$$ \partial_{\nu,\mu}\left( e^{\pi i\,  \trace(A Z)}\right) = 
   (2\pi i) A_{\nu,\mu}\,  e^{\pi i\,  \trace(A Z)}
$$
with $A = (A_{\nu,\mu})$. 
\qed
 
The proof of Theorem \ref{Siegelthm} will follow directly from 
the following Proposition and Theorem \ref{Ibuthm}. 
\begin{prop}
\label{thprop}
Let $F$ and $F'$ be Siegel modular forms of degree $2$ and weight 
$k$ and $k'$, respectively. Assume that $F$ and $F'$ can be written 
as theta series with harmonic coefficients.
Then
$[F,F']_l$ (defined in Theorem \ref{Siegelthm}) is a Siegel modular 
form of degree $2$ and  weight $k+k'+2l$.
\end{prop}
{\it Proof.}
Let $F$ and $F'$ be Siegel modular forms of degree $n=2$ and 
weight $k$ and $k'$, respectively which can be written as 
theta series with spherical coefficients.
Then one has 
$$
F  = \sum_{G\in M_{2m,n}(\Z)}  Q (S^ {1/2}G )\, 
     e^{\pi i \, \trace(G^t S  G  Z)}, \qquad
F' = \sum_{G'\in M_{2m',n}(\Z)} Q'(S'^{1/2}G')\,
     e^{\pi i \, \trace(G'^tS' G' Z)} 
$$
for some symmetric, positive, even and unimodular matrices 
$S$ and $S'$, and spherical polynomials $Q$ and $Q'$ of weight 
$d$ and $d'$, respectively and $k=m+d$ and $k'=m'+d'$.

Then, from the very definition of $[\cdot,\cdot]_l$ and 
Lemma \ref{trilem}, we have 
$$
[F,F']_l = (2\pi i)^{n l} \sum_{G,G'} 
              e^{\pi i\, \trace((G^t S G+G'^t S' G')Z)} \
              \tilde P(S^{1/2}G,S'^{1/2}G') 
$$
with 
$\tilde P(X,X')= Q(X)\, Q'(X')\, P(X,X')$
where
$$P(X,X') = \sum_{r+s+p=l} 
            C_{r,s,v}(k,k') 
            \det(X^t X)^r 
            \det(X'^t X')^s 
            \det(X^t X + X'^t X)^p.
$$ 
We will now show that the polynomial $\tilde P$ is spherical
of weight $w=d+d'+2l$ in the matrix variable 
$$ Y = \left( \begin{array}{cc}  X\\ 
                                 X'
              \end{array}\right).$$
If this is the case Theorem \ref{thetathm} implies that $[F,F']_l$ 
is a Siegel modular form of degree $2$ and weight $k+k'+2l$.

Note that $\tilde P$ clearly satisfies the first property 
in the definition of a spherical polynomial with $w=d+d'+2l$.
Therefore,  we only have to show that  
$$ \Delta_Y \tilde P = (\Delta_{S^{1/2} G} + \Delta_{S'^{1/2}G'}) 
    \tilde P = 0.$$

To calculate $\Delta_Y \tilde P$ we use a change of basis such that 
$S$ and $S'$ become equal to $\id_{2m}$ and $\id_{2m'}$, respectively. 
To show that $\tilde P$ is spherical we need to know the following 
expressions which can easily be obtained 
by a straightforward calculation
(here we use that $n=2$ for the first time)
\begin{eqnarray}
\Delta_G( \det(G^t G) ) &=&  4(m-1/2)\, \trace(G^t G),
 \nonumber \\
(\grad_G( \det(G^t G) ))^2 &=& 4\,\trace(G^t G)\, \det(G^t G),
 \nonumber \\
\Delta_G( \det(G^t G+G'^t G')) &=& 4(m-1/2)\, \trace(G^t G)+4k\,\trace(G'^tG'),
 \nonumber \\
(\grad_G( \det(G^t G+G'^t G')))^2 &=& 
   4\, \trace(G^t G+G'^t G')\, (\det(G^t G) - \det(G'^t G')) 
 \nonumber \\
 &&   + 4\, \trace(G'^t G')\, \det(G^t G+G'^t G'),
 \nonumber \\
\grad_G( \det(G^t G) ) \cdot \grad_G( \det(G^t G+G'^t G') &=& 
    4\, \trace(G^t G)\, \det(G^t G),
\nonumber\\
\grad_G(Q(G)) \cdot \grad_G( \det(G^t G) ) &= &
 2d\, Q(G)\, \trace(G^t G),
\nonumber\\
\grad_G(Q(G)) \cdot \grad_G( \det(G^t G+G'^t G'))&= &
2d\, Q(G)\, \trace(G^t G+G'^t G').
\nonumber
\end{eqnarray}
Note the for deriving the last two expressions one has to use the 
second property in the definition of spherical polynomials
with 
$$
A=\left( \begin{array}{cc} \lambda &0\\ 0&\lambda \end{array}\right),\quad
A=\left( \begin{array}{cc} 1&0\\ \lambda & 1\end{array}\right)
\quad {\rm{and}}\quad
A=\left( \begin{array}{cc} 1&\lambda\\ 0 & 1\end{array}\right).
$$
Let us, for the moment, use 
$D=\det(G^t G)$, $D' =\det(G'^t G')$ and 
$D'' =  \det(G^t G+G'^t G')$.
Then the last equalities imply that 
\begin{eqnarray}
\frac{1}{4}(\Delta_G + \Delta_{G'})&
                 \Big( Q(G) Q'(G') D^r  D'^s  D''^p \Big)\qquad \quad= &
\nonumber\\
  &r(\alpha +r)\, \trace(G^t G) &  
                           Q(G) Q'(G') D^{r-1}  D'^s   D''^p
\nonumber\\
+ &s(\beta  +s)\, \trace(H^t H) & 
                           Q(G) Q'(G') D^r D'^{s-1}   D''^p
\nonumber\\
+ &p(\gamma+l+p)\, \trace(G^t G+H^t H) &  
                           Q(G) Q'(G') D^r D'^s   D''^{p-1}.
\nonumber   
\end{eqnarray}
Here we have used $\alpha = k-3/2$, $\beta = k'-3/2$ and 
$\gamma = k+k'-3/2$.
With these expressions it is clear that the equation 
$\Delta_X \tilde P=0$ is certainly satisfied if the coefficients 
$C_{r,s,p}$ obey
\begin{eqnarray}
0 &=& (r+1)(\alpha+r+1)\, C_{r+1,s,p}(k,k') + 
      (p+1)(\gamma+l+p+1)\, C_{r,s,p+1}(k,k') 
\nonumber\\
0 &=& (s+1)(\beta +s+1)\, C_{r,s+1,p}(k,k') + 
      (p+1)(\gamma+l+p+1)\, C_{r,s,p+1}(k,k'). 
\nonumber
\end{eqnarray}
It is simple to check that these conditions are indeed satisfied by 
the $C_{r,s,p}$ given in the theorem.
\qed
\begin{rem}
Note that a covariant bilinear differential operator which can be
written in terms of the operator $\D$ is, by the recursion relations 
obtained  at the end of the proof of Proposition \ref{thprop},
equal to a multiple of $[\cdot,\cdot]_l$.
Furthermore, the recursion relations imply the explicit form of 
the combinatorial factors $C_{r,s,p}$.
\end{rem}
\smallskip\noindent
{\it Proof of Theorem \ref{Siegelthm}.}
Firstly, note that the proof of Proposition \ref{thprop} 
implies for $Q=Q'=1$ that the polynomial $\tilde P = P$ is spherical 
of weight $2l$. Secondly, it can be written as 
$$P(X,X') =  \tilde Q(X^t X, X'^t X')$$
with
$$\tilde Q(a,b) = \sum_{r+s+p=l} 
                  C_{r,s,v}(k,k') 
                  \det(a+b)^p
                  \det(a)^r 
                  \det(b)^s 
$$ 
so that it satisfies the assumptions of Theorem \ref{Ibuthm}.
Finally, note that $ [F,F']_l$ can be written as 
$$ [F,F']_l =  {\cal D}(F(Z) F(Z'))|_{Z=Z'} $$ 
where 
${\cal D} = \tilde Q(\partial_{\nu,\mu},\partial_{\nu',\mu'}')$.
This implies the desired result.
\qed

\newpage
\section{Concluding remarks and open questions}
\label{Conclusion}

Let us end with some remarks and mention some open questions.

\begin{enumerate}
\item  It is obvious from the proofs of the Theorems \ref{Jacobithm}
       and \ref{Siegelthm} that they also hold true for the case of 
       Jacobi forms on  $\Gamma\semi\Z^2\subset \SL\semi\Z^2$ and
       Siegel forms on $\Gamma'\subset\SG2$ if $\Gamma$ and $\Gamma'$ 
       are finite index subgroups of $\SL$ and $\SG2$, respectively.
\item  In the generic case, {\it i.e.}\ where $k$ and $k'$ are large enough,
       the dimension of the space of Rankin-Cohen operators for Jacobi
       forms is given by B\"ocherer's result so that Theorem \ref{Jacobithm}
       describes a basis of this space explicitly.
       This implies that, for $k$ and $k'$ large enough, the dimension of 
       the space of Rankin-Cohen operators from 
       $S^{(2)}_k\times S^{(2)}_{k'}$ to $S^{(2)}_{k+k'+v}$
       is one dimensional for even $v$ and zero otherwise.
       (This can be verified by noting that any Rankin-Cohen operator for
       Siegel modular forms induces a Rankin-Cohen operator for Jacobi
       forms via the Jacobi-Fourier expansion (Theorem \ref{Siegelexp}). 
       Hence it is enough to show that there is only one (up to 
       multiplication by a constant) Rankin-Cohen operator for Jacobi 
       forms for even $v$ and none 
       for odd $v$ which can  be `lifted' to a Rankin-Cohen operator for 
       Siegel modular forms. 
       In the general case, this can indeed be done using Theorem 
       \ref{Jacobithm}.)
\item  Using the relation between Jacobi forms and modular forms 
       of half-integral weight one can show that the operators 
       $(\frac{d}{dX})^{v/2}[\cdot,\cdot]_{X,v}$ ($v\in 2\N$)
       can be obtained from the Rankin-Cohen operators for elliptic 
       modular forms (for more details see \cite{Ch1}). It seems that 
       this is the only Rankin-Cohen operator for Jacobi forms for which 
       such a result holds true.
\item  The operators $(\frac{d}{dX})^{v/2}[\cdot,\cdot]_{X,v}$ ($v\in 2\N$)
       can be used to define generalised Rankin-Cohen algebras \cite{CE} 
       which have very similar properties to the Rankin-Cohen algebras in 
       the elliptic case considered in ref.\ \cite{Z}.
\item
       It would be interesting to understand how our constructions
       (via generating functions or theta series)
       can be generalised to higher Jacobi and  Siegel modular
       forms and multilinear differential operators.
       We hope to discuss this in a future publication.
\item  Is it possible to obtain the explicit formulae for the
       Rankin-Cohen operators in the case of Siegel modular froms 
       from the representation theory of $\mbox{Sp}(4,\R)$?
       In this context the Rankin-Cohen operators can be viewed
       as certain projection operators. \item
       Is it possible to obtain the dimension of the space of 
       covariant bilinear operators from $S^{(2)}_k\times S^{(2)}_{k'}$ 
       to $S^{(2)}_{k+k'+v}$ not using Jacobi forms?
       Of course one possibility would be to use 
       theta series and rephrase the question in terms of 
       the theory of invariants.
\item
       Is there any connection between covariant bilinear operators 
       for Jacobi forms and/or Siegel modular forms and automorphic 
       pseudodifferential operators like in the case of elliptic 
       modular forms considered in ref.\ \cite{CMZ}?
\end{enumerate}

\newpage
{\bf Acknowledgements}

We would like to thank R.E.\ Borcherds and D.\ Zagier for 
valuable discussions and S.\ B\"ocherer and T.\ Ibukiyama 
for making their results available prior to publication.

This work was done during the first author was visiting the 
Department of Pure Mathematics and Mathematical Statistics 
at the University of Cambridge. Y.C. would like to thank the 
Department and, in particular, J.\ Coates for support and 
warm hospitality during her stay in Cambridge.



\begin{thebibliography}{[20]}
\bibitem{Boe}
S.\ B\"ocherer, private communication.

\bibitem{Ch1}
Y.\ Choie, 
{\it Jacobi Forms and the Heat Operator},
preprint (1996), Math. Z. (to appear).

\bibitem{Ch2}
Y.\ Choie, 
{\it Jacobi Forms and the Heat Operator II},
preprint (1996).

\bibitem{CE}
Y.\ Choie, W.\ Eholzer,
{\it Generalised Rankin-Cohen Algebras and Jacobi Forms},
in preparation.

\bibitem{C}
H.\ Cohen, 
{\it Sums Involving the Values at Negative Integers of $L$ Functions
     of Quadratic Characters},
Math. Ann. 217 (1975), 271-285.

\bibitem{CMZ}
P.B.\ Cohen, Y.\ Manin, D.\ Zagier,
{\it Automorphic Pseudodifferential Operators},
preprint MPI-95-132 (1995).

\bibitem{EZ}
M.\ Eichler, D.\ Zagier,
{\it The Theory of Jacobi Forms},
Prog. Math. Vol. 55,
Birkh{\"a}user, Boston-Basel-Stuttgart (1985).

\bibitem{Frei}
E.\ Freitag,
{\it Siegelsche Modulfunktionen},
Grundlehren der mathematischen Wissenschaften 254,
Springer, Berlin--Heidelberg--New York (1983).

\bibitem{Ibu}
T.\ Ibukiyama,
{\it On Differential Operators on Automorphic Forms
     and Pluri-Harmonic Polynomials},
preprint 1991 (unpublished).

\bibitem{R}
R.\ Rankin,
{\it The Construction of Automorphic Forms from the 
     Derivatives of a given Form},
J. Indian Math. Soc. 20 (1956), 103-116.

\bibitem{Z}
D.\ Zagier,
{\it Modular Forms and Differential Operators}, 
Proc. Indian Acad. Sci. (Math. Sci.) Vol. 104, No. 1 (1994), 57-75.


\end{thebibliography}
\end{document}